\documentclass[]{article}
\begin{document}
\newtheorem{proposition}{Proposition}[section]
\newtheorem{definition}{Definition}[section]
\newtheorem{lemma}{Lemma}[section]

\title{\bf Some Results about Triangular Representations of Lie Algebras}
\author{Keqin Liu\\Department of Mathematics\\The University of British Columbia\\Vancouver, BC\\
Canada, V6T 1Z2}
\date{June , 2014}
\maketitle

\begin{abstract} We introduce the concept of a triangular  representation of a Lie algebra, give a counterpart of Ado's theorem, and discuss $2$-irreducible triangular modules over a nonreductive 
Lie algebra. \end{abstract}

\bigskip
It is well-known that the representation theory of semisimple Lie algebras is a fruitful area in Lie 
theory, but  it is basically not known how to develop the finite dimensional representation theory of non semisimple Lie algebras. This 
fact makes us feel strongly that the ordinary definition of a representation of a Lie algebra is too 
general to give rise 
to an effective way of developing the representation theory of non semisimple Lie algebras. If a 
restrictive definition 
has to be used to study representations for non semisimple Lie algebras, which conditions should be added 
to the ordinary definition of a representation of a Lie algebra? 

\medskip
By a fact which is due to E. Cartan, a nonreductive Lie algebra does not have a faithful finite dimensional irreducible representation. Hence, replacing the linear transformations in the general linear Lie 
algebra by the linear 
transformations which have more than two invariant subspaces should be one change we should do to introduce a restrictive definition for a representation of a non semisimple Lie algebra. By Ado's theorem, a finite dimensional Lie algebra over 
a field of characteristic $0$ has a finite dimensional faithful representation such that every element in its nil 
radical ( which is the largest nilpotent ideal) can be represented by a nilpotent linear transformation. Since it is 
easy to get a lot of nilpotent linear transformations on a finite dimensional graded vector space and a finite 
dimensional module over a finite dimensional Lie algebra automatically has a graded vector space structure after treating the module as a module over a Levi factor of the Lie algebra, replacing the general vector space in the ordinary definition of a representation of a Lie algebra by a graded vector space and requiring that the elements in the nil radical of a Lie algebra are represented in terms of a class of obvious nilpotent linear transformations on the graded vector space should be another change we should do. Based on these ideals, we write this paper to introduce the 
concept of a triangular  representation of a Lie algebra and to give some evidences which indicate that it is 
reasonable to develop the representation theory of non semisimple Lie algebras by using triangular  representations.

\medskip
This paper consists of two parts. In part 1, we use linear transformations which have many invariant 
subspaces to define the general triangular linear Lie algebra. It turns out that the general triangular 
linear Lie algebra,  which is a subalgebra of the general linear Lie algebra, has a natural graded vector space structure. Based on this property of the general triangular linear Lie algebra, we add two more conditions in the ordinary definition of a representation of a Lie algebra to introduce the concept of 
a triangular  representation of a Lie algebra. We finish Part 1 with the main result of this paper which claims that every finite dimensional Lie algebra over 
a field of characteristic $0$ has a finite dimensional faithful triangular representation. Because of this counterpart of Ado's theorem, we know that although the concept of a triangular  representation of a Lie algebra is restrictive, this concept is still general enough. This is our reason to 
initiate the study of the triangular representations of Lie algebras. In Part 2, we discuss finite 
dimensional $2$-irreducible triangular representations of a nonreductive Lie algebra 
$s\ell_2^{\Lambda}$, 
where  $s\ell_2^{\Lambda}$ is the semidirect product of the $3$-dimensional Lie algebra 
$s\ell_2$ by the radical 
$rad(s\ell_2^{\Lambda})$ of $s\ell_2^{\Lambda}$ such that $rad(s\ell_2^{\Lambda})$ is both an Abelian subalgebra of 
$s\ell_2^{\Lambda}$ and a finite dimensional irreducible $s\ell_2$-module.

\medskip
Throughout this paper, $\mathbf{F}$ is a field of characteristic $0$, all vector spaces are finite dimensional vector spaces over $\mathbf{F}$, and the notation $k\in \mathbf{Z}_{\ge j}$ means that $k$ is an integer which is larger than or equal to $j$.

\bigskip
\bigskip
\noindent
\setcounter{section}{0}
{\bf \large{Part 1: The counterpart of Ado's theorem}}

\bigskip
For convenience, we use $\bigoplus_{k\ge j}V_k$ to denote a graded vector space 
$\bigoplus_{k\in \mathbf{Z}_{\ge j}}V_k$. If $V=\bigoplus_{k\ge j}V_k$ is a finite dimensional graded 
vector space, then 
$V_k=\{0\}$ for large $k$.

\begin{definition}\label{def1.1} Let $V=\bigoplus_{k\ge 0}V_k$ be a finite dimensional graded vector space. A {\bf triangular linear transformation} $f: V\to V$ is a linear transformation on $V$ satisfying the condition: $f\big(\bigoplus_{k\ge j}V_k\big)\subseteq \bigoplus_{k\ge j}V_k$ for $j\in \mathbf{Z}_{\ge 0}$.
\end{definition}

Let $End^{\bigtriangleup}\big(\bigoplus_{k\ge 0}V_k\big)$ be the set of all triangular linear transformations on the graded vector space $\bigoplus_{k\ge 0}V_k$. We write 
$g\ell^{\bigtriangleup}\big(\bigoplus_{k\ge 0}V_k\big)$ for 
$End^{\bigtriangleup}\big(\bigoplus_{k\ge 0}V_k\big)$ viewed as a subalgebra of the general linear Lie 
algebra $g\ell\big(\bigoplus_{k\ge 0}V_k\big)$ and call it the {\bf general triangular linear Lie algebra}.

\begin{proposition}\label{pr1.1} If $V=\bigoplus_{k\ge 0}V_k$ is a finite dimensional graded vector space,
then $g\ell^{\bigtriangleup}(V)=\bigoplus_{j\ge 0}g\ell^{\bigtriangleup}_j(V)$ as the direct sum of vector spaces, where $g\ell^{\bigtriangleup}_j(V)$ with $j\in \mathbf{Z}_{\ge 0}$ is defined by
$$
g\ell^{\bigtriangleup}_j(V):=\{\, f\,|\, \mbox{$f\in g\ell^{\bigtriangleup}(V)$ and 
$f(V_k)\subseteq V_{k+j}$ for $k\in \mathbf{Z}_{\ge 0}$}\,\}.
$$
\end{proposition}

Obviously, every element in $\bigoplus_{j\ge 1}g\ell^{\bigtriangleup}_j(V)$ is a nilpotent linear transformation. Using these nilpotent linear transformations in 
$\bigoplus_{j\ge 1}g\ell^{\bigtriangleup}_j(V)$, we now introduce the concept of a triangular  representation of a Lie algebra.

\begin{definition}\label{def1.2} Let $\mathcal{L}$ be a Lie algebra and let $V=\bigoplus_{k\ge 0}V_k$ be a finite dimensional graded vector space. A {\bf triangular  representation} of the Lie algebra $\mathcal{L}$ 
is a Lie algebra homomorphism $\varphi: \mathcal{L}\to g\ell^{\bigtriangleup}(V)$ which satisfies the 
following two more conditions:
\begin{description}
\item[(i)] $\varphi(\mathcal{L}_{semi})\subseteq g\ell^{\bigtriangleup}_0(V)$, where $\mathcal{L}_{semi}$
is a Levi factor of $\mathcal{L}$;
\item[(ii)] $\varphi(nil\,\mathcal{L})\subseteq \bigoplus_{j\ge 1}g\ell^{\bigtriangleup}_j(V)$, where 
$nil\,\mathcal{L}$ is the nil radical of $\mathcal{L}$.
\end{description}
\end{definition}

The condition (i) in Definition \ref{def1.2} makes it possible to use the fruitful representation theory of semisimple Lie algebras in the study of triangular  representations of non semisimple Lie algebras. Because
the condition (i) in Definition \ref{def1.2} depends on a Levi factor of a Lie algebra, a triangular  representation $\varphi: \mathcal{L}\to g\ell^{\bigtriangleup}(V)$ is also called the triangular  representation of $\mathcal{L}$ with respect to the Levi factor $\mathcal{L}_{semi}$, and the graded vector space $V$ is called the {\bf triangular $\mathcal{L}$-module} with respect to the Levi factor 
$\mathcal{L}_{semi}$. By the following proposition, the concept of triangular  representation of a
Lie algebra does not depends on the choice of a Levi factor of the Lie algebra.

\begin{proposition}\label{pr1.2} Let $\mathcal{L}$ be a finite dimensional Lie algebra. If a finite 
dimensional graded vector space $V=\bigoplus_{k\ge 0}V_k$ is a triangular $\mathcal{L}$-module with 
respect to a Levi factor of the Lie algebra $\mathcal{L}$, then $V$ is also a triangular 
$\mathcal{L}$-module with respect to any Levi factor of the Lie algebra $\mathcal{L}$.
\end{proposition}

The next proposition proves that the concept of triangular  representation of a Lie algebra is in fact quite natural.

\begin{proposition}\label{pr1.3} A finite dimensional Lie algebra can be made into a graded vector space such that its adjoint representation is a triangular  representation.
\end{proposition}

We finish Part 1 of this paper with the following counterpart of Ado's theorem.

\begin{proposition}\label{pr1.4} A finite dimensional Lie algebra over 
a field of characteristic $0$ has a finite dimensional faithful triangular representation.
\end{proposition}

\bigskip
\bigskip
\noindent
{\bf \large{Part 2:  Triangular $2$-irreducible $s\ell_2^{\Lambda}$}-modules}

\bigskip
Let $\mathcal{L}$ be a finite dimensional Lie algebra and $n\in \mathbf{Z}_{\ge 0}$. A finite dimensional graded vector space $V=\displaystyle\bigoplus_{k=0}^n V_k$ with respect to a Levi factor 
$\mathcal{L}_{semi}$ of $\mathcal{L}$ is said to be {\bf $(n+1)$-irreducible} if each $V_k$ is an irreducible 
$\mathcal{L}_{semi}$-module for $k=0$, $1$, ... , $n$. 

\medskip
The purpose of Part 2 of this paper is to describe $2$-irreducible triangular $s\ell_2^{\Lambda}$-modules,
where the nonreductive Lie algebra $s\ell_2^{\Lambda}$ with $\Lambda\in \mathbf{Z}_{\ge 1}$ is defined as follows:
$$
s\ell_2^{\Lambda}:=\mathbf{F}f\oplus \mathbf{F}h\oplus\mathbf{F}e\oplus
\mathbf{F}z_0\oplus\mathbf{F}z_1\oplus\cdots \oplus\mathbf{F}z_{\Lambda},
$$
$$
[h, \, e]=2e,\qquad [h, \, f]=-2f,\qquad [e, \, f]=h,\qquad [z_j, \, z_{j'}]=0,
$$
$$
[h, \, z_j]=(\Lambda-2j)z_j,\qquad [f, \, z_j]=z_{j+1},\qquad [e, \, z_j]=j(\Lambda-j+1)z_{j-1},
$$
where $\Lambda\ge j,\,j'\ge 0$ and $z_j:=0$ for $j\not\in\{0,\, 1,\, ..., \Lambda\}$.

\medskip
Clearly, the $3$-dimensional Lie algebra $s\ell_2=\mathbf{F}f\oplus \mathbf{F}h\oplus\mathbf{F}e$ is a Levi factor of $s\ell_2^{\Lambda}$, and $\mathbf{F}z_0\oplus\mathbf{F}z_1\oplus\cdots \oplus\mathbf{F}z_{\Lambda}$
is both the radical and the nil radical of $s\ell_2^{\Lambda}$.

\medskip
For arbitrary four nonnegative integers $m$, $n$, $s$ and $N$ which satisfy the conditions:
\begin{equation}
m+2s=\Lambda +n+2N,\qquad n\ge s\ge 0, \qquad m\ge N\ge 0,
\end{equation}
and arbitrary $n-s$ scalars $a_1$, ... , $a_{n-s}\in\mathbf{F}$, there is a  $2$-irreducible triangular 
$s\ell_2^{\Lambda}$-module $\mathcal{M}_{m,\, n}^{s,\, N}(a_1, ... , a_{n-s})$ which is defined as follows:
$$
\mathcal{M}_{m,\, n}^{s,\, N}(a_1, ... , a_{n-s})=
\Big(\displaystyle\bigoplus_{i=0}^n \mathbf{F} u_i\Big)\bigoplus
\Big(\displaystyle\bigoplus_{k=0}^m \mathbf{F} w_k\Big);
$$
$$
h\cdot u_i=(n-2i)u_i,\quad f\cdot u_i=u_{i+1},\quad e\cdot u_i=i(n-i+1)u_{i-1}\quad
\mbox{for $n\ge i\ge 0$};
$$
$$
h\cdot w_k=(m-2k)w_k,\quad f\cdot w_k=w_{k+1},\quad e\cdot w_k=k(n-k+1)w_{k-1}\quad
\mbox{for $m\ge k\ge 0$};
$$
$$
z_j\cdot w_k=0 \quad \mbox{for $\Lambda\ge j\ge 0$ and $m\ge k\ge 0$};
$$
$$
z_j\cdot v_i=0 \quad \mbox{for $\Lambda\ge j\ge 0$,  $n\ge i\ge 0$ and $s-1\ge i+j$};
$$
\begin{eqnarray*}
&&z_j\cdot v_{s-j+\theta}=\left(\sum_{k=0}^{\theta}(-1)^{j-k}{j \choose k}\,
\frac{(m-N-\theta+k)!}{(N+\theta -k)!\,m!}\,a_{_{\theta -k}}\right)w_{_{\theta+N}}\\
&&\mbox{for $min\{\, s+\theta,\, \lambda\,\}\ge j\ge 0$ and $n-s\ge \theta\ge 0$};
\end{eqnarray*}
\begin{eqnarray*}
z_j\cdot v_{n-j+\theta}=\left(\sum_{k=0}^{j-\theta}(-1)^{j-\theta -k}{j \choose \theta +k}
\frac{(m-N-n+s+k)!}{(N+n-s -k)!\,m!}a_{n-s -k}\right)w_{_{n-s+\theta +N}}&&\\
\mbox{for $\lambda\,\ge j\ge\theta \ge 1$},
\qquad\qquad\qquad\qquad\qquad\qquad\qquad\qquad\qquad\qquad\qquad\qquad\qquad&&
\end{eqnarray*}
where $a_0:=1$, $\displaystyle{j \choose k}:=\displaystyle\frac{j!}{k!(j-k)!}$ for $j\ge k$, 
$\displaystyle{j \choose k}:=0$ for $k> j$, $u_i:=0$ for $i\not\in\{0,\, 1,\, ..., n\}$, and
$w_k:=0$ for $k\not\in\{0,\, 1,\, ..., m\}$.
 
\medskip
We finish this paper with the following

\begin{proposition} Let $\mathcal{M}$ be a $2$-irreducible triangular 
module over the non-reductive Lie algebra $s\ell_2^{\Lambda}$. If 
$rad(s\ell_2^{\Lambda})\cdot \mathcal{M}\ne\{0\}$, then 
$\mathcal{M}=\mathcal{M}_{m,\, n}^{s,\, N}(a_1, ... , a_{n-s})$
for some  $n-s$ scalars $a_1$, ... , $a_{n-s}\in\mathbf{F}$ and some non-negative integers 
$m$, $n$, $s$ and $N$ which satisfy the conditions in (1).
\end{proposition}

\medskip
{\bf Acknowledgment.} I thank Professor Robert V. Moody for his help and support.

\end{document}